\documentclass[a4paper,12pt]{article}
\usepackage{amsmath,amssymb,amsfonts,amscd,amsthm}
\usepackage{graphicx}
\usepackage{dsfont}
\usepackage{blindtext}

\usepackage{accents}

\usepackage[titletoc, toc, title, page]{appendix}
\usepackage{color}
\usepackage{enumerate}
\usepackage{fancyhdr}
\usepackage{parskip}
\usepackage{lastpage}
\usepackage[utf8]{inputenc}
\usepackage[top=2cm, bottom=2cm, left=2cm, right=2cm]{geometry}
\usepackage{subcaption}
\usepackage{bbm}
\usepackage{pgf,tikz}
\usepackage{lipsum}
\usepackage{authblk}
\usetikzlibrary{arrows}
\usepackage{color}

\usepackage{xcolor}

\definecolor{bluegray}{rgb}{0.4, 0.6, 0.8}

\definecolor{dgreen}{rgb}{0.1,0.6,0.1}
\definecolor{bluegreen}{rgb}{0.1,0.5,0.2}

\def\bg{\color{bleu1}}
\def\bk{\color{black}}

\usepackage{xcolor}

\newtheorem{theorem}{\bg Theorem}                             

\newtheorem{proposition}{\bg Proposition}

\newcommand{\A}{\mathbb A}
\newcommand{\C}{\mathbb C}
\newcommand{\R}{\mathbb R}

\newcommand{\Z}{\mathbb Z}
\newcommand{\N}{\mathbb N}

\newcommand{\T}{\mathbb{T}}

\renewcommand{\min}{\mathop{\mbox{min}}}



\renewcommand{\th}{\theta} 
 
\newcommand{\om}{\omega} 
\newcommand{\cO}{\mathcal O} 

\numberwithin{equation}{section}

\title{\bg \large{Convergence of the Birkhoff normal form sometimes  \\
implies 
convergence of
a normalizing transformation }}
\date{}
\author[1]{\small{Rafael de la Llave}\thanks{Supported in part by NSF    DMS 1800241}}
\author[2]{\small{Maria Saprykina}\thanks{Supported in part by  the Swedish
    Research Council, VR 2015-04012}}
\affil[1]{School of Mathematics, Georgia Institute of Technology}
\affil[2]{Dept. of Mathematics, KTH Royal Institute of Technoligy}

\newcommand{\Addresses}{{
  \bigskip
  \footnotesize

  \textit{E-mail addresses:} 

  R.~de la Llave: \texttt{rafael.delallave@math.gatech.edu}

 M.~Saprykina: \texttt{masha@kth.se}

}}

\usepackage{sectsty,ulem} \definecolor{bleu1}{RGB}{0,57,128}
\subsectionfont{\color{bleu1}}
\sectionfont{\color{bleu1}}
\subsectionfont{\color{bleu1}}

\usepackage[utf8]{inputenc}
\usepackage[english]{babel}

\usepackage[colorlinks,linkcolor=bleu1,anchorcolor=green,citecolor=bleu1]{hyperref}

\usepackage{amsmath,amsthm,amssymb,mathrsfs,latexsym,amsfonts}
\usepackage{bbm}
\usepackage{amssymb}
\usepackage{amsmath}

\usepackage{color}

\theoremstyle{plain}
\newtheorem{thm}{Theorem}

\newtheorem{Rem}[thm]{Remark}
\newtheorem{Lem}[thm]{\bg Lemma}

\def\D{{\mathbb D}}
\def\Z{{\mathbb Z}}
\def\T{{\mathbb T}}
\def\R{{\mathbb R}}
\def\C{{\mathbb C}}
\def\N{{\mathbb N}}

\newcommand{\dt}{{\delta}}
\newcommand{\eps}{{\epsilon}}
\newcommand{\kp}{{\kappa}}
\newcommand{\lb}{{\lambda}}

\newcommand{\si}{{\sigma}}
\def\th{{\theta}}
\newcommand{\Om}{{\Omega}}

\newcommand{\cL}{{\mathcal L}}

\def\<{\langle}
\def\>{\rangle}
\newcommand\beq{\begin{eqnarray}}
\newcommand\eneq{\end{eqnarray}}
\newcommand\beal{\begin{aligned}}
\newcommand\enal{\end{aligned}}
\newcommand\blm{\begin{Lem}}
\newcommand\elm{\end{Lem}}
\newcommand\brm{\begin{Rem}}
\newcommand\erm{\end{Rem}}
\usepackage{calc}
\usepackage{accents}

\newcommand{\comm}[1]{}
\newcommand{\comment}[1]{}

\usepackage{color}

\definecolor{orange}{rgb}{1,0.5,0}

\newcommand{\vertiii}[1]{{\left\vert\kern-0.25ex\left\vert\kern-0.25ex\left\vert #1 
    \right\vert\kern-0.25ex\right\vert\kern-0.25ex\right\vert}}

\begin{document}

\maketitle

\begin{abstract}
Consider an analytic Hamiltonian system near its analytic invariant
torus $\mathcal T_0$ carrying zero frequency. 
We assume that the Birkhoff normal form of the Hamiltonian at
$\mathcal T_0$ is convergent and 
has a particular form: it is an analytic function of its
non-degenerate quadratic 
part. We prove that in this case there is an analytic canonical 
transformation---not just a formal power series---bringing 
the Hamiltonian into its Birkhoff normal form.
 


\end{abstract}
\maketitle

\section{Introduction} 

The goal of this paper is to study the convergence of 
the transformations of an analytic Hamiltonian system in a 
neighborhood of an invariant torus to the Birkhoff normal form. 
Here we assume that
the frequency vector at the invariant torus is very resonant, hence already at the formal level, the existence of the 
Birkhoff normal form has obstructions.  
The main result, Theorem~\ref{main} below, will show that if the obstructions
for the formal equivalence between the system and its Birkhoff normal form vanish and the normal form is convergent and has
a particular form, 
then the system is analytically equivalent to its normal form.  Hence, this result 
can be considered as a part of the rigidity program: identifying 
obstructions for a weak form of equivalence whose vanishing 
implies a stronger form of equivalence.

\subsection{Classical theory of normal forms: existence and uniqueness.} 
Consider an analytic function   
\beq\label{Ham1}
H(I,\theta)=\left< \lb_0,I\right> +\cO^2(I), 
\eneq
where $\theta\in \T^d=\R^d/\Z^d$, $I\in (\R^d,0)$, $\left< \cdot,\cdot\right>$
denotes the usual scalar product in $\R^d$, and $\lambda_0\in \R^d$ is a constant
vector called the {\it frequency vector}. The Hamiltonian system
associated to it is
$
\dot I=\partial_{\theta} H(I,\theta), \
\dot \theta =- \partial_{I}H(I,\theta)
$.
Note that we are assuming the standard symplectic form. In particular, the set $\mathcal T_0:=\{0\}\times \T^d
$ is an invariant torus of this system.
 We say that $H(I,\theta)$ {\it has a Birkhoff normal form  (BNF)} $N(I)$ in
a neighborhood of $\mathcal T_0$ if $N(I)$ is a  {\it formal} power series, and there
exists a {\it formal} symplectic transformation $\Psi(I,\theta)$,
tangent to the identity
$$
\Psi(I,\theta)=(I+\cO^2(I),\theta+\cO(I))
$$
such that 
$$
H\circ \Phi(I,\theta)= N(I) 
$$
in the sense of formal  power series. 
Any canonical coordinate change $\Phi(I,\theta)$  as above is called a {\it normalizing
transformation}. The following 
fundamental result is called the Birkhoff normal form \cite{SiegelM71, MeyerHO}. 
For $H(I,\theta)$ as above, assume that $\lambda_0$ 
satisfies a Diophantine condition: there
exist constants $(C,\tau)$ such that for all 
$k\in \Z^d\setminus \{0\}$ we have
\beq\label{DC}
\left| \<\lb_0, k \> \right| \geq
C|k|^{-\tau}.
\eneq
Then $H(I,\theta)$ has a (formal) Birkhoff normal form. Moreover, if a
normal form exists and $\lambda_0$ is rationally independent, then 
the Birkhoff normal form is unique (up to trivial changes relabelling the 
actions). Note that the normalizing transformations are not unique, since 
composing $\Phi(I,\theta)$ with any transformation that preserves $I$ gives a normalizing transformation.

Birkhoff normal form is an important tool in the study of Hamiltonian
systems. Already the assumption of existence and nondegeneracy of the normal
form has strong dynamical consequences (see, e.g., \cite{EFK15} Th.C).
The importance of the BNF becomes even stronger if the normal form is convergent, and even
more so if there exists an analytic normalizing transformation.  

The standard way of constructing BNF, which we will review in more detail 
later, is to proceed iteratively, devising transformations 
that normalize $H(I,\theta)$ up to  the coefficients of order $I^n$. 
The normalization step involves solving 
differential equations with analytic conditions.  The 
Diophantine conditions \eqref{DC} can be somewhat 
weakened to subexponential growth ($ \lim_{N \to \infty}\frac{1}{N} 
\log  \sup_{|k| \le N } \left| \<\lb_0, k \> \right|^{-1} = 0$).

If  $\lb_0$ is resonant, one cannot guarantee the existence of the 
Birkhoff normal form even at the level of 
formal power series, since 
there may be  some terms in the formal power series of $H$ that cannot 
be eliminated by a canonical transformation. On the other hand, 
there are, of course, systems (e.g the BNF itself, or changes of variables 
from it) for which one can construct a BNF even in the resonant case. 
Then one speaks of the 
Birkhoff-Gustavson normal form \cite{Gustavson66}.


Analogous definitions and statements hold true for symplectic maps
in a neighborhood of a fixed point. Even if the formal elimination procedures
are very similar, the analysis is very different. 
Handy references for the classical theory of Birkhoff normal forms 
are \cite{SiegelM71, MeyerHO,Murdock,EFK13, EFK15}.




\subsection{Generic divergence both of the Birkhoff 
Normal Form and the normalizing transformation.} 

The BNF and the normalizing transformations are constructed as 
formal power series. 
The following natural questions are of great importance: the first one is whether the BNF 
converges for Hamiltonians in a certain class. 
The second---whether there is 
a convergent normalizing transformation. 

Concerning the first question, R. Perez-Marco \cite{PM} proved the following dichotomy:   for any given nonresonant quadratic
part, either the BNF is generically divergent
or it always converges. The original proof was done in the setting of 
Hamiltonian systems having
a non resonant elliptic fixed point.  The extension of this result to the case of the torus, that  
is not completely straightforward,  has been worked out 
by  R. Krikorian, see Theorem 1.1 in \cite{Kri}.

Up to very recently it was unclear which of the possibilities is actually realized. A large progress has 
been made  by R. Krikorian \cite{Kri}, who proved that there exists a real analytic symplectic diffeomorphism $f$ of a two-dimensional annulus
such $f(\T \times \{0\})=(\T \times \{0\})$, $f(\theta,0)=(\theta+\om_0,0)$ with $\om_0$ Diophantine and having a 
non-degenerate {\it divergent} Birkhoff normal form. Combined with the aforementioned result of Perez-Marco,  
this implies that Birkhoff Normal Form of an analytic Hamiltonian is 	 ``in general" divergent.

Concerning the normalizing
transformations, H. Poincar\'e proved that they are divergent for a 
generic Hamiltonian. 
C.~L.~Siegel proved the same statement in a neighborhood of an
elliptic fixed point (in fact, for a larger class of Hamiltonians than
just generic, \cite{Siegel54}).
This is implied by showing that the orbit structure of the map 
in any neighborhood is 
very different from that of the Birkhoff normal 
form (which is integrable). Analogous results for 
symplectic maps near an elliptic fixed point appear in \cite{Russmann59}. 
Very different arguments showing divergence of normalizing transformations 
for generic systems appear in \cite{Zehnder73} and for 
some concrete polynomial mappings in \cite{Moser60}.

\subsection{Convergence of the transformations 
under the Diophantine conditions for some particularly simple BNF}

There are classes of Hamiltonians for which we can guarantee the convergence of the normalizing transformation. 
The following influential  rigidity  result  was proved independently by A.~D.~Bruno \cite{Br71} and 
H.~R\"ussmann \cite{Russmann67}. Note that the main assumption is 
that the (in principle only formal) BNF is of  a particular kind. 

 
Consider an analytic Hamiltonian $H(I,\theta)$, whose
frequency $\lb_0$ satisfies a Diophantine condition \eqref{DC}.
 Assume moreover that the Birkhoff normal form $N(I)$ of $H(I,\theta)$ is a  {\it formal} function $B$ of one single variable  
$\Lambda_0:=\< \lb , I\>$, i.e., 
$$
N(I)=B(\Lambda_0(I)).
$$ 
Then there exists an {\it analytic} normalizing 
transformation,  and the BNF is, in fact, analytic.

We remark that Bruno proves the above result under a weaker condition on $\lb_0$ than \eqref{DC}. 
For analogous statements in the case of invariant tori see \cite{Br89}. Other modifications can be found in
\cite{Russmann02, Russmann04}. 
This result has been recently generalised to a much more general
context by Eliasson, Fayad and
Krikorian \cite{EFK13, EFK15}.
We stress that in all these works mentioned above, $\lb_0$ is assumed to be non-zero and
the crucial assumption is that $\lb_0$ satisfies a Diophantine-type condition
and that the BNF is of a very simple form.

\subsection{``Sometimes'' convergence of the BNF implies 
convergence of a normalizing transformation.}

Our main result is close in spirit to the above works, 
but it {\it does not rely on a Diophantine condition}. In fact, 
we consider a special class of diffeomorphisms 
such that the frequency $\lb_0$ is zero. Thus, 
the BNF is  degenerate in the
previous sense. But within this class of Hamiltonians we 
just use a standard non-degeneracy
assumption on the quadratic part. Namely, we prove the following.

\bigskip

\begin{theorem}\label{main}

Assume the following:
\begin{itemize}
\item[$(A_1)$] $H(I,\theta)$ has a formal Birkhoff normal 
form $N(I)$ that starts with
quadratic terms in $I$, i.e
there exists a {\it formal} symplectic change of variables 
$\Psi(I, \th)$, tangent to the identity, i.e. 
$
\Psi(I,\theta)=(I+\cO^2(I),\phi+\cO(I))
$, such that 
$$
H\circ \Psi (I,\th) = N(I)=N_0(I)+ \cO^3(I)
$$
in the sense of power series.

\item[$(A_2)$] $N_0(I)=I^{tr} \Omega  I$ (for some symmetric $\Omega$) is 
non-degenerate: $\det \Om \neq 0$.

\item[$(A_3)$] $N(I)=B(N_0(I))=N_0 + \sum_{j=2}^\infty b_j (N_0(I))^j$ where $B$ is an 
analytic function.
\end{itemize}

Then there exists an invertible {\it analytic} symplectic transformation
$$
\Phi(I,\theta)=(I+\cO^2(I),\phi+\cO(I))
$$  
such that
\begin{equation}\label{conjugacy}
H\circ \Phi (I,\th) = N(I).
\end{equation}
\end{theorem}

\bigskip

 Note that we start from a resonant torus, so that the existence of a BNF of the form we assume, requires vanishing of (formal) obstructions.
Hence, our main result can be reformulated as saying that the formal assumptions imply convergence of the normalizing transformation.
 
Similar rigidity statements have appeared in other contexts. In 
\cite[Ch. 5]{Poincare92}, H.Poincar\'e studied the 
formal power series of canonical 
transformations, which send  a family of Hamiltonian 
systems into a family of integrable systems (in the sense of 
power series). In  \cite{Poincare92} 
it was shown  that these formal 
power series do not exist unless there are some conditions
(which
are not met in the three body problem for arbitrary masses).
 The non-existence of 
formal power series, a fortiori implies the non-existence of 
analytic families of analytic transformations integrating the 
three body problem.

The paper \cite{Ll} proved a converse to the result in 
\cite{Poincare92}: if the system satisfies {a  very specific 
and generic}  non-degeneracy condition, 
then, existence of a formal power series that integrates the family
of transformations in the sense of power series implies existence of a convergent one. 

Assumption $A_3$ is there for technical purposes, see Sec.~\ref{s_simplif}. Note that it is trivial for $d=1$. This assumption reminds of that
of R\"ussmann in \cite{Russmann67, Russmann02, Russmann04}.

The assumption that the Birkhoff normal 
form is a function of $N_0$ 
has been discussed in \cite{Ga} under the name of 
{\it relative integrability}.  Two Hamiltonian dynamical 
systems are 
relatively integrable when one of them can be  obtained  
from the other by a symplectic change of coordinates and 
a reparameterization of the time which only depends on 
the total energy. That is, the orbit structures of the two systems in an 
energy surface are equivalent up to a change of scale of 
time. The paper \cite{Ga} includes several arguments for why 
the notion of relative integrability is natural when discussing 
formal equivalence. 
In the present paper, however, the focus lies on the notion of
 equivalence under a  symplectic change of variables. 
We show that, for a certain class of systems, equivalence in 
the sense of formal power series implies equivalence 
in the sense of analytic canonical changes of variables. 
Hence, our main result can be understood as a rigidity result.  
The class of systems for which this rigidity 
result holds can be succinctly described as the set of systems that 
relatively integrable with respect to the main term.

\medskip
In the context of formal equivalence implying analytically convergent 
equivalence, it is natural to formulate: 

\noindent{\bf Conjecture.  }{\it  Assume that an analytic Hamiltonian $H(I,\theta)$  as in \eqref{Ham1}
has a convergent  BNF  that satisfies the  non-degeneracy assumption
that the frequency map is a local diffeomormphism.
Then there is 
a convergent normalizing transformation. }

\medskip

Note that the problems studied in \cite{Russmann67} and \cite{Br71} do not satisfy 
the hypothesis of the conjecture, even though they satisfy the conclusions. 

In the other direction one can construct examples \cite{S} of 
analytic maps near a hyperbolic
fixed point such that the Birkhoff normal form is quadratic (in the
above notations, $N=\Lambda_0$) with a
non-resonant set of eigenvalues, and any normalizing transformation to
the normal form diverges. In these examples, 
the eigenvalues form carefully chosen Liouville vectors.
That is, the paper \cite{S} shows that, depending on the 
Diophantine conditions, quadratic normal forms may be rigid or not. 
The models in \cite{S} do not satisfy the hypothesis of 
the  conjecture above.


\subsection{Overview of the proof.} 
The standard method of obtaining the Birkhoff Normal form is  
an iterative procedure in which we construct the 
transformations order by order: at the $n$-th step of the procedure 
one computes the $n$-th order terms in the Taylor expansions,
assuming that all the terms of lower orders are computed. 
It would appear natural to follow this scheme and 
try to {\it estimate} the transformations at each step of the 
recursive procedure. Unfortunately, this seems technically unfeasible. 
One of the main complications in any possible  proof of 
convergence of the transformations  is that even if 
the BNF is unique, the formal transformations $\Phi_N$ are very far from 
unique (Since the BNF depeds only on the actions, the $\Phi_N$
can be composed with  any canonical 
transformation which moves the angles but preserves the actions. So, 
an essential ingredient of any proof of convergence should be 
a especification of how to choose the normalizing transformations.

In this paper we use a 
quadratically convergent method in which we double 
the number of known coefficients at each step.  
Roughly -- see more details in 
the next paragraphs --  we will show that if the formal obstructions 
vanish we can choose a sequence of canonical transformations 
that proceed to converge quadratically: doubling the order of 
the BNF at every step of the construction. More importantly, 
there is a specific choice of the transformation that satisfies 
very explicit bounds. The bounds on the new transformation in terms 
of the remainder turn out to involve a loss of derivatives. Therefore
we need to implement a Nash-Moser scheme to estimating the important objects in a sequence of 
domains which decrease slowly.

Here is a short overview of the proof; all the necessary notations are introduced in the next 
section. At the $n$-th step of the iterative procedure we will start with a
Hamiltonian of the form
$$
H_n(I,\th)=N_n(I)+ \widetilde{R_n}(I,\th),
$$
where $N_n(I)$ is a polynomial in $I$ of degree $m_n=2^n+1$, and the remainder term $\widetilde{R_n}$ is small in the following sense: 
for a certain domain-dependent norm, introduced in Sec. \ref{sec_norms}, for a certain small $\dt_n$ (we assume $\dt_n \to 0$
with $n\to \infty$) and $\kp>0$ 
the remainder term satisfies
$|\widetilde{R_n}|_{\rho_n,\rho_n}\leq \delta_{n}^\kp $.

At this step we construct a symplectic change of coordinates $\Phi_n$,
such that 
$$
H_n\circ\Phi_n (I,\th)=N_{n+1}(I)+ \widetilde{R_{n+1}}(I,\th),
$$
where $N_{n+1}$ has degree $m_{n+1}=2m_n-1$, and 
$|\widetilde{R_{n+1}}|_{\rho_{n+1},\rho_{n+1}}\leq \delta_{n+1}^\kp
=2^{-\kp}\delta_{n}^\kp$. 

We construct $\Phi_n$ as a time one map of a the flow of
a Hamiltonian vector field $F_n$.
The main ingredient consists in constructing and estimating the norm
of $F_n$ (and thus $\Phi_n$), which is found as a solution of a
certain homological equation (see \eqref{eq_homol} and in a
simplified form \eqref{homol-eq}). In general,
this equation may not have even a formal solution
unless some constraints are met. However, the assumption of Theorem
\ref{main} implies that this equation does have a formal solution.
The key observation in this paper is the following:
if this homological equation has a {\it formal solution}, then it also
has an {\it analytic solution with  tame
estimates} for it
(in the sense of Nash-Moser theory). This statement is the contents of 
Lemma \ref{lemma-homol-eq0}. 
We note that the tame estimates use an argument different from the 
matching of powers.

The procedure can be repeated, because the main assumption used 
to show the existence of solutions of the Newton equation is 
that there is a formal solution to all orders. This assumption 
is clearly preserved if we make any analytic change of variables. 
Once we know that the Newton procedure can be repeated infinitely 
often, the convergence is more or less standard.

\section{Notations and a step of induction. }

\subsection{ Notations.}


\subsubsection{Norms and majorants.} \label{sec_norms} 

Let $\T^d=\R^d /\Z^d$ be a $d$-dimensional torus, and for $\si>0$ consider its complex extension
$\T^d_\si =\left(\R^d+(-\si,\si)\sqrt{-1} \right)/ \Z^d$. Let $\D^d_{\rho }=\{I\in
\C^{d}: |I|< \rho  \}$ be a complex disc, and define the "d-dimensional annulus"
$$
\A_{\rho,\si }:=\D_{\rho}^d\times \T^d_\si .
$$
Let $\cO(\A_{\rho,\si })$ be the set of functions holomorphic in $\A_{\rho,\si }$ that are real symmetric, i.e., such that
$\overline{f(\bar I,\bar \th)}=f(I, \th)$ (where the bar stands for the complex conjugate).
We use supremum norms over $\A_{\rho,\si }$, denoted by $\|f\|_{\rho,\si 
}$. In the same way we define the set $\cO(\D_{\rho })$ 
with the corresponding norm 
$\|f\|_{\rho }$ being the sup-norms over the disc
$\D^d_{\rho }$.

For a function $f\in \cO(\A_{\rho,\si })$ consider its Taylor-Fourier representation in the powers of $I$:
$
f(I,\th)=\sum_{j\in \N^d} \sum_{k\in \Z^d} f_{j,k}e^{2\pi i
  \<k,\theta \>} I^j
$.
Consider a majorant for $f$ of the form 
$$
\widehat f(I) = \sum_{j\in \N^d} \sum_{k\in \Z^d} |f_{j,k}| I^j e^{2\pi|k|\si}
$$ 
We denote 
by $|f|_{\rho,\si }$ the norm of the corresponding majorant $\widehat f(I)$:
$$
|f|_{\rho,\si }=\|\widehat  f\|_{\rho,\si }.
$$
Clearly, $\|f\|_{\rho,\si } \leq |f|_{\rho,\si }$. Analogous notation
$|f|_{\rho}$ corresponds to the norm $\|f\|_{\rho}$ above.

In what follows we will mostly have $\si =\rho$.

\subsubsection{Important constants for the iterative procedure.}\label{not_const}

\begin{itemize}
\item Let $\rho_0=\min \{1, \rho\}$,
 

\item The order of polynomials involved in the $n$-th step of the iterative procedure is
$$
m_{n} =2^n +1.
$$

\item The norm of the rest term $\widetilde{R_n}$ at the $n$-th step will be
estimated as $|\widetilde{R_n}|_{\rho_n}\leq \delta_{n}^\kp $.
Let 
$$
\begin{aligned}
&\kappa = d + 6, \\
&b =  2^{-(\kappa + 3)} \\
& \delta_0= \rho_0  b 2^{-3} = \rho_0  2^{-(\kappa + 6)}\\
&  \delta_{n+1}= 2^{-1} \delta_{n}.
\end{aligned}
$$ 

\item  Finally, let 
$$
q_n=(2b)^{2^{-(n+1)}},
$$
and
$$
\rho_{n+1}=(\rho_n-3\dt_n)q_n.
$$

\end{itemize}
\subsubsection{Polynomials.}\label{not_poly}
In the iterative procedure we will work with polynomials in
$I$ whose coefficients depend on $\th$. 
\begin{itemize}

\item Let
\begin{equation}\label{def_N0}
N_0(I)=I^{tr} \Omega  I
\end{equation}
where $\Omega$ is 
a symmetric non-degenerate matrix: $\det \Om \neq 0$.

\item An expression $M=f(\th) I^k $ (where $k$ is a multi-index) is called a monomial.

\item We will say that a monomial $M_{k,l}=I^ke^{2\pi i \< l,\th\>}$ is {\it resonant} if it satisfies $\{N_0, M\}=0$.

\item $R^{[j]} (I,\th)$ stands for a homogeneous polynomial in $I$ of
degree $j$ with coefficients depending on $\th$:  
$$
R^{[j]}(I,\th)=\sum_{|k |=j} r_{k }(\th ) I^{k }.
$$

\item We also use notation $R^{[m,n]}$ to denote  the range of degrees in $I$:
$$
R^{[m,n]} (I,\th)=\sum_{j=m}^n R^{[j]} (I,\th), \quad R^{[\geq m]} (I,\th)=\sum_{j=m}^\infty R^{[j]} (I,\th).
$$
\end{itemize}

 Let $m_n$ be as above. The following functions will be of special importance. 

\begin{itemize}

\item The normal form  $N(I)$ is assumed to have the form
\begin{equation}\label{N=B}
N(I) = B(N_0(I))= N_0(I) + \sum_{j=2}^\infty b_j (N_0(I))^j.
\end{equation} 
Denote
\beq\label{def_N}
N_n=N^{[2,m_n]}= \left( B(N_0)\right)^{[2,m_n]};
\eneq
in particular, since $m_0=2$, $N_0=N_0^{[2,m_0]}=N_0^{[2]}$ is quadratic.

\item The rest term at the $n$-th inductive step is $\widetilde{R_n}(I, \theta)$:
\beq\label{def_tR}
\widetilde{R_n}= \widetilde{R_n}^{ [>m_{n} ] }.
\eneq

\item We will also need polynomials in $I$ with $\theta$-dependent coefficients:
$R_n(I, \theta)$ and $F_n(I, \theta)$ of the following degrees:
\beq
R_n= R_n^{ [m_n+1,m_{n+1} ] }, \quad  F_n=F_n^{[m_n, m_{n+1}-1]}.
\eneq

\end{itemize}

\subsection{Base of induction: an equivalent problem.}
\begin{Lem}\label{lemma_base}
Suppose that 
$$
H(I,\th)=N_0(I)+ \widetilde {R_0} (I,\theta) \in
\cO(\A_{\rho,\si }),
$$ 
where $|\widetilde{R_0}|_{\rho,\si }\leq \delta$, and  there exists a
formal (resp., analytic) symplectic transformation
$$
\Psi(I,\theta)=\left(\phi(I,\theta),\,
  \psi(I,\theta) \right)=(I+\cO^2(I),\theta + \cO(I))
$$
such that 
$$
H\circ\Psi(I,\th)= N(I)=N_0(I) + \sum_{j=2}^\infty b_j (N_0(I))^j .
$$
 
Then for any $a > 0$ there exists a  Hamiltonian $\widehat
H(I,\th)$ and a formal (resp., analytic) symplectic transformation
$
\widehat  \Psi(I,\theta)=(I+\cO^2(I),\theta + \cO(I))
$ such that
$$
\widehat H\circ \widehat  \Psi(I,\th) =N_0(I)+ \widehat   {R_0} (I,\theta) \in
\cO(\A_{\frac1{a} \rho,\si }),
$$
where $| \widehat {R_0}|_{\frac1{a} \rho,\si }\leq a \delta$, and  
$$
N(I)=N_0(I) + \sum_{j=2}^\infty b_j a^{2(j-1)}(N_0(I))^j .
$$
 \end{Lem}

\bg \noindent {\it Proof. \ }\bk Define $\widehat
H(I,\th)=\frac1{a^2}H(aI,\th)$, and  
$
\widehat  \Psi(I,\theta)= \left( \frac{1}{a}\phi(aI,\theta),\,
  \psi(aI,\theta) \right)
$.
It can be verified directly that $\widehat  \Psi$
is symplectic and tangent to the identity. Moreover,
$$
\widehat  H \circ \widehat  \Psi(I,\theta)= \frac{1}{a^2}
H(\phi(aI,\theta), 
\psi(aI,\theta)) = 
N_0(I) + \sum_{j=2}^\infty b_j a^{2(j-1)}(N_0(I))^j.
$$
\hfill {\bg $\Box$}

\bigskip

\subsection{Induction step.} 
While the base of induction is given by formula \eqref{est_R0},
the step of the iterative procedure is provided by the following
proposition. 

\begin{proposition}\label{lemma1}
For a fixed $n > 0$, let $m_n$, $\rho_n$ and $\dt_n$ be as in Sec. \ref{not_const} above.
Suppose that $H_n(I, \theta)$ is formally conjugated to the BNF of the
form \eqref{N=B}:
$$
N(I)=N_0(I) + \sum_{j=2}^{\infty} b_j (N_0(I))^j,
$$
and the normal form satisfies:
\begin{equation}\label{est_N}
|N^{[ m_n+j ]}|_{\rho_n}<\dt_n^{\kappa+1}, \quad j=0, \dots m_n,
\end{equation}
and denoting
$
g_{2j}(I) = jb_j (N_0(I))^{j-1}
$,
we assume
\begin{equation}\label{est_g}
|g_j|_{\rho_n} \leq \frac1{4^j}, \quad j=1,\dots ,m_n ;
\end{equation}

Suppose that
$$
H_n(I, \theta) = N_n(I) + \widetilde{R_n} (I, \theta),
$$
where $N_n(I) = \left( B(N_0(I))\right)^{[2,m_n]} $
and  $\widetilde{R_n}=  \widetilde{R_n}^{ [>m_{n} ] } $ satisfies 
$$
|\widetilde {R_n}  |_{\rho_n,\rho_n} \leq \delta_n^{\kappa}.
$$

Then there exists a symplectic change of coordinates $\Phi_n:(I',
\theta')\mapsto (I,\theta)$,
$$
\Phi_n(I', \theta')=(U^{(n)} (I', \theta'), V^{(n)}(I', \theta')),
$$
given by a Hamiltonian $F_n=F_n^{[m_{n},m_{n+1}-1]}$ such that 
\beq\label{eqHn}
H_{n+1}(I', \theta'):=H_n \circ \Phi_n(I', \theta') = N_{n+1} (I') +\widetilde { R_{n+1} }(I', \theta'),
\eneq
where $N_{n+1}(I')=N^{[2,m_{n+1}]}(I')$, $\widetilde { R_{n+1}}(I',\theta')=\widetilde {R_{n+1}}^{[>m_{n+1}]}(I',\theta')$, and
\beq\label{estRn}
|\widetilde { R_{n+1}} |_{\rho_{n+1},\rho_{n+1}} \leq \dt_{n+1}^{\kappa} . 
\eneq

Moreover, $\Phi_n(I', \theta')=(U^{(n)}(I', \theta'), V^{(n)}(I',
\th'))$ satisfies 
\beq\label{estPhi}
\sum_{j=1}^d \|U_j^{(n)}(I', \theta')-I'_j \|_{\rho_n-3\dt_n,\rho_n-3\dt_n} + 
\|V_j^{(n)}(I', \theta')-\theta'_j \|_{\rho_n-3\dt_n,\rho_n-3\dt_n} < \dt_n,
\eneq
and the inverse map, 
$\Phi_n^{-1}(I, \theta):=(U^{(-n)}(I, \theta),V^{(-n)}(I, \theta))$, satisfies
\beq\label{estPhi-1}
\sum_{j=1}^d \| U_j^{(-n)}(I, \theta)-I_j \|_{\rho_n-3\dt_n,\rho_n-3\dt_n} + 
\| V_j^{(-n)}(I, \theta)-\theta_j \|_{\rho_n-3\dt_n,\rho_n-3\dt_n} < \dt_n.
\eneq
\end{proposition}

The proof of this proposition constitutes the
main technical tool of this paper. It
implies Theorem 1 in a standard way. See, e.g., \cite{Russmann67},
pp 61-63). For convenience, we give a proof below. 

\subsection{Proof of Theorem 1.} 
Lemma \ref{lemma_base} permits us to assume without loss of 
generality that for the given Hamiltonian 
$
H_0(I,\th):=H(I,\th) = N_0(I)+\widetilde{R_0}(I,\th),
$
\begin{equation}\label{est_R0}
|\widetilde{R_0}|_{\rho_0,\rho_0}\leq \delta_{0}^\kp. 
\end{equation}
Since  the function $B$ is analytic,  the same lemma permits us to assume that  \eqref{est_N} and \eqref{est_g} hold for each $n$.

The step of induction is provided by Proposition \ref{lemma1}.
Since $H_{n}$ is formally reducible to the normal form $N$, the same can be
said about $H_{n+1}$. 

Repetition of this process leads to a sequence of transformations 
$$
T_n=\Phi_0\circ\Phi_1\circ\dots\circ\Phi_{n-1}.
$$ 
Let us show that $T_n$ converges to
the desired coordinate change $\Phi=T_\infty$, analytic in the polydisc
$\A_{\rho_\infty,\rho_\infty}$, where $\rho_0 b < \rho_\infty < \rho_0$. Indeed, with the notations of Sec. \ref{not_const},
$$
3 \sum_{k=0}^\infty \dt_k \leq 3\cdot 2\dt_0 < 3\cdot 2 \rho_0 b
2^{-3} <\rho_0 b .
$$
Then 
for any $n$ we have 
$$
\rho_{n+1}=q_n (\rho_{n} - 3\dt) \geq  \rho_{0} \prod_{j=0}^n q_j - 3\sum_{j=0}^n \dt_n
\geq \rho_{0} \prod_{j=0}^\infty q_j - 3\sum_{j=0}^\infty \dt_n\geq \rho_{0} 2b - 3\cdot 2\dt_0 >b\rho_{0}
$$
It is left to prove that $T_n$ converges of to an analytic function
$T_\infty$, satisfying \eqref{conjugacy}. 
Denote the variables, involved in the $n$-th step of the induction by 
$w_{n-1}=(I,\theta)$ and $w_{n}=(I',\theta')$, where
$$
w_{n}=\Phi_{n-1}^{-1} w_{n-1}.
$$
In these notations,
$$
w_{0}=\Phi_{0} \circ \Phi_1 \circ \dots  \circ\Phi_{n-1}w_n = T_n w_n.
$$
Now, for $w_{n}=(I',\theta')$ we have 
$$
H\circ T_n(I',\theta')= N_n(I')+\widetilde{R_n}(I',\theta').
$$
Since $(\Phi_{n}(I',\theta')-(I',\theta'))$ 
starts with the terms of degree $2^n$ in $I'$, for each $j$ the expansion of
 $(T_n(I',\theta')-T_{n+j}(I',\theta'))$ starts with the terms of
 degree $2^n$ in $I'$. 
This implies that
the sequence of maps $T_n$ formally converges, when $n\to \infty$, to
a formal map $T_\infty$ such that \eqref{conjugacy} holds:
$$
H\circ T_\infty (I',\theta')= N(I').
$$
We still need to show that $T_\infty$ is analytic. It is more convenient to prove that the maps
$$
T_n^{-1}:=\Phi_{n-1}^{-1}\circ\dots\circ\Phi_1^{-1}\circ\Phi_{0}^{-1}
$$
converge to an analytic map $T_\infty^{-1}$. 

By Proposition \ref{lemma1}, the map
$$
w_{n+1}=\Phi_n^{-1} w_n
$$
is analytic  in $\A_{\rho_0 b/2,\rho_0 b/2}$, and  for all $n$ we have:
$$
|\Phi_n^{-1} w_n- w_n|_{\rho_0 b/2,\rho_0 b/2}\leq \dt_n,
$$
since $\rho_n-3\dt\geq \rho_{n+1}> \rho_0b$ for all $n$.
Therefore, the map $T_n^{-1}$ such that
$$
w_{n}=T_n^{-1} w_0 
$$
is analytic in $\A_{\rho_0 b/4,\rho_0 b/4}$, and for such $w_0$ we have
$$
|T_n^{-1} w_0|\leq \sum_{j=0}^{n-1} |T_{j}^{-1}(w_j)-w_j| +|w_0|\leq 
\sum_{j=0}^{\infty}\dt_j +\rho_0b/4 \leq \rho_0b/2.
$$
Estimate
$$
|T_{n+m}^{-1}(w_0)-T_n^{-1}(w_0)|_{\rho_0 b/4,\rho_0b/4}\leq 
\sum_{j=n}^{n+m-1} |T_{j}^{-1}(w_j)- w_j)|_{\rho_0 b/4,\rho_0b/4}\leq 
\sum_{j=n}^{\infty}\dt_j =2^{1-n}\dt_0
$$
implies the convergence of the sequence of maps $T_n^{-1}$ to an
analytic map $T_\infty^{-1}$ in $\A_{\rho_0 b/4,\rho_0 b/4}$.
Since  the formal inverse of $T_\infty^{-1}$ is the series $T_{\infty}$, the latter also defines an analytic function,  
providing the desired coordinate change. We set $\Phi=T_\infty$ in the notations of Theorem \ref{main}.
\hfill {\bg $\Box$}


\section{Formal analysis.}
Here we start the proof of Proposition \ref{lemma1} by the formal analysis of the
iterative procedure.

\subsection{Iterative Procedure.}

Given $H_n$ as in  Proposition \ref{lemma1},
we will construct $\Phi_n$ as the time one map of the flow of a
Hamiltonian $F_n$, i.e., $\Phi_n = X_{F_n}^1$ where  $X_{F_n}^t$ is
the flow
defined by 
$$
\dot I =F_\th(I,\th),\quad \dot \th =-F_I (I,\th).
$$
In this case, $\Phi_n$ is automatically symplectic.

Notice that the normalising transformation $\Phi_n$, as well as the
corresponding generating function $F_n$, is not unique (one can
compose with rotations in the angles which preserve the  actions, for
example). Clearly, the transformation that converges has to be
very carefully chosen.

In the following Lemma~\ref{lem_help} we show
that if a (formal) normalizing transformation exists, then there
exists (another) normalizing transformation of a special kind. Namely, such that the corresponding generating function 
is a polynomial (in the sense of section \ref{not_poly}),
$F_n=F_n^{[m_n,     m_{n+1}-1]}$ and free from resonant monomials (see notations in Sec. \ref{not_poly}).

The idea of the proof is that we can always move the formal normalizing
transformation  by composing with some  transformations that
do not change the normal form. Therefore, we can ensure that the
normalizing transformations belong to a space which is transversal
to the space spanned by resonant monomials.
Note that in the proof of Lemma~\ref{lem_help}
we use crucially the fact that the normal form is
a function of $N_0$ so that the resonant terms are the same at all orders.

There are some analogies between Lemma~\ref{lem_help}  and 
Proposition 2.6 in \cite{Ll}, but  that result is significantly less
delicate since there is an extra parameter that controls the smallness. 
In our case, the variable $I$  controls both  the smallness and the distance
to the origin at the same time.

Let $\{\cdot, \cdot\}$ denote the standard 
Poisson bracket. Recall that 
for a differentiable function $G$ it holds:
$$
\frac{d}{dt}G\circ X_F^t =\{ G,F\}\circ X_F^t.
$$

\bigskip

\begin{Lem}\label{lem_help}
Suppose that  for $H(I,\th)$ 
there exist $N_{2m}(I)=N_0 + B(N_0)$ with $B(X)=\sum_{j=2}^{m} b_j X^{j}$,  $R(I,\th)=R^{[> 2m]}(I,\th)$ and $G(I,\th)=\mathcal O^2(I) $ such that
$\Psi:=X_G^1$ satisfies 
\[
H\circ \Psi (I,\th)=N_{2m}(I)+R(I,\th).
\]

\def\tR{{\tilde R}}
\newcommand{\dbtilde}[1]{\accentset{\approx}{#1}}
\def\ttR{{\dbtilde {R}}}
\def\tPsi{{\tilde \Psi}}
\def\tG{{\tilde G}}
\def\cL{{\mathcal L}}

\begin{enumerate}

\item Then there exists $\tG(I,\th)$, that is free from resonant monomials of order $< 2m$, 
  such that $\tilde \Psi :=X_{\tilde G}^1$ normalises $H$ to the same normal form, i.e.,
  for some 
 $\tR(I,\th)=(\tilde R)^{[> 2m]}(I,\th)$ we have:
$$
H\circ \tPsi(I,\th)=N_{2m}(I)+\tR (I,\th).
$$

\item If, an addition to the previous assumption, we have that the original $H(I,\th)$ has the form 
$$
H(I,\th)= N_{m}(I)+ R^{[>m]}(I,\th),
$$
where $N_{m}=N_{m}^{[2,\dots , m]}$, then there exists a polynomial $F=F^{[m, 2m-2]}$, that is free from resonant monomials, 
 such that $\Phi :=X_{F}^1$ normalises $H$ to the same normal form, i.e., for some 
 $\ttR(I,\th)=\ttR^{[> 2m]}(I,\th)$ we have:
$$
H\circ \Phi (I,\th)=N_{2m}(I)+\ttR(I,\th).
$$

\end{enumerate}
\end{Lem}

\bg \noindent{\it Proof:   (1).} \bk      
All the calculations below are made in the sense of formal Taylor-Fourier expressions.
Suppose that $K(I,\th)$ is such that $\{N_0, K\}=0$. Notice that in this case
$\{N_{2m}, K\}= B' (N_0)\{N_0, K\} =0$.
Use $K(I,\th)$ as a  Hamiltonian to define $k(I,\th):=X_{K}^1$. 
Then by Taylor formula  we have:
$$
\begin{aligned}
H\circ \Psi \circ k & =  (N_{2m}+R )\circ k= (N_{2m} + R )\circ
X_{K}^t \Big|_{t=1} = N_{2m} + R + \{(N_{2m}+R), K\}  \\
&+ \frac12\{\{(N_{2m}+R), K \} ,K\}+ \dots =
N_{2m}+R_1,
\end{aligned}
$$
where $R_1(I,\th)=R_1^{[> 2m]}(I,\th)$.   


It is a classical fact that the composition $\Psi \circ k$
in the sense of formal power series is the time-one
map of another Hamiltonian given by the Cambell-Baker-Dynkin
formula \cite[Appendix C]{Dragt}, \cite[Appendix]{LlMM}; here we denote it by CBD formula.
Note that in these references the usual  notation for the Hamiltonian vector field defined by $G$ is $\cL_G$, and $\exp(\cL_G)$ stands for its time one map. 
In the present paper the same map is denoted by $X_G^1$. Now, suppose that
$\Psi = X_G^1$ and $ k = X_K^1$. 
CBD formula implies that the composition of these maps satisfies:
\[
\begin{split} 
& \tilde \Psi :=\Psi \circ k  = X_{\tilde G}^1 , \quad \text{where}\\ 
& {\tilde G} =  G + K + \frac{1}{2}\{G, K\}
+ \frac{1}{12}\{G, \{G, K\}\}
- \frac{1}{12}\{K, \{K, G\}\} + \cdots 
\end{split}
\]
The last sum is to be understood in the sense of formal power series in $I$.

To prove Lemma \ref{lem_help}, we use CBD formula, and choose $K$ recursively (order by order in $I$) so that ${\tilde G} $
has no resonant terms up to order $2m$. At each step of the recursion we choose $(-K(I,\theta))$ 
to be equal to the lowest order resonant term of $G$, and set ${\tilde G} $ to be the new $G$. As we saw above, the map $\tilde \Psi =\Psi \circ K$, used as a normalization map,   brings $H$ to the same normal form as $\Psi$ did. But its generating Hamiltonian ${\tilde G} $ has no lower order resonant monomials. Iterating this procedure, we get a normalization with the desired property.

\medskip

\bg \noindent{\it (2).} \bk 
Since we can normalise $H=N_{m}+R^{[>m]} $ to $N_{2m}$ with the help of the generating function 
  $G=\mathcal O^2(I)$, then, by {\it (1)}, we can also achieve the 
normalization using the transformation $\tilde \Psi$ generated by a resonance-free Hamiltonian
  $\tilde G $.
  Note that $\tilde G =\mathcal O^2(I)$.
  
By the Taylor formula for power series,   we have:
\[
\begin{aligned}
H\circ \tilde \Psi   & =  (N_{m}+R^{[>m]} ) \circ {\tilde \Psi} = (N_{m} + R^{[>m]} )\circ
X_{\tilde G}^t \Big|_{t=1} = N_{m} + R^{[>m]} \\
&+ \{(N_{m}+R^{[>m]}), \tilde G\}  
+  \frac12\{ \{(N_{m}+R^{[>m]}), \tilde G \}, \tilde G \}+ \dots  
=
N_{2m}+R_1.
\end{aligned}
\]

Since $\tilde G $ is resonance-free, any monomial $P$ in  $\tilde G $ gives a non-zero impact $\{N_0, P\} $ to the sum above, whose order in $I$ is strictly larger than the order of $P$.
By comparing the orders of the coefficients in $I$ we see that the lowest possible order of a monomial in 
$\{N_{0}, \tilde G\} $ is the same as that in $R^{[>m]}$, and hence $\tilde G=\tilde G^{[\geq m]}$. 
Finally notice that the reduced 
generating function $F:=\tilde G^{[m, 2m-2]}$ produces the same normal form.

\hfill {\bg $\Box$}

\bigskip

The following lemma introduces the notations used in the proof of the Main Theorem. Here we use the results of Lemma \ref{lem_help} to relate the conjugating function to the
solutions of the homological
equation \eqref{eq_homol} below.

 \medskip

\begin{Lem}\label{lem_homol_eq}
Adopt the notations for the degrees of polynomials from 
Sec. \ref{not_poly} (in particular,  $N_n=N^{[2,m_n]}$ as in \ref{def_N}, and $R_n=R_n^{ [m_n+1,m_{n+1} ] }$).
Let $B(X)=\sum_{j=1}^{\infty} b_j X^{j}$. Suppose that $H_n$ has the form 
$$
H_n= N_n +\widetilde{ R_n} = N_n + R_n +\widetilde{ R_n}^{[>m_{n+1}]}.
$$
where $N_{n}(I)=N_0 + B(N_0)^{[4,m_n]}$.

Suppose that there exists $G(I,\th)=\mathcal O^2(I) $ such that
$\Psi:=X_G^1$ satisfies 
\[
H\circ \Psi (I,\th)=N_{m+1}(I)+R(I,\th).
\]

Then there exists a polynomial (in $I$)
$F_n=F_n^{ [m_n,m_{n+1}-1] }$ with the following properties: 
 the time one map $\Phi_n:= X_{F_n}^1$ satisfies 
$$
H_{n+1}:=H_n \circ \Phi_n  = N_{n+1} + \widetilde{ R_{n+1}},
$$
and $F_n$ satisfies
\begin{equation}\label{eq_homol}
\{N_n, F_n\}^{ [m_n+1,m_{n+1}] }+ R_n + N_n - N_{n+1} =0 ,
\end{equation}
and 
$$
\widetilde{ R_{n+1}}:= A_n+B_n +C_n,
$$
where
\begin{equation}\label{not_AnBn}
A_n:= \widetilde{R_n}^{[>m_{n+1}]}\circ \Phi_n  , \quad 
B_n:=\int_0^1\{ (1-t) \{N_n, F_n \}+R_n,F_n\}\circ X_{F_n}^t dt,
\end{equation}
\begin{equation}\label{not_Cn}
C_{n} = (\{N_n, F_n  \} )^{[>m_{n+1}]} .
\end{equation}
\end{Lem}
Notice that  the expressions for $A_n$, $B_n$, $C_n$ 
start with terms of order $m_{n+1}+1$, and hence, $\widetilde{ R_{n+1}}=\widetilde{ R_{n+1}}^{[>m_{n+1}]}$, as needed.

\medskip 

\bg \noindent{\it Proof.} \bk  
Let $m=m_n=2^n+1$. Then $m_{n+1}= 2m-1$.  With the notations for the degrees of polynomials from 
Sec. \ref{not_poly},
Lemma \ref{lem_help} implies that  there exists a polynomial $F_n=F_n^{ [m_n,m_{n+1}-1] }$ 
such that $\Phi_n:= X_{F_n}^1$ satisfies $H_n \circ \Phi_n = N_{n+1} +\widetilde{R_{n+1}}$.
By the Taylor formula we have:
\beq\label{eq1}
\begin{aligned}
H_n \circ \Phi_n  = &(N_n + R_n +\widetilde{ R_n}^{[>m_{n+1}]})\circ
X_{F_n}^t \Big|_{t=1} = N_n + \{N_n, F_n\} + R_n +\\ 
&\int_0^1\{ (1-t)
\{N_n, F_n \}+R_n,F_n\}\circ X_{F_n}^t \, dt + \widetilde{R_n}^{[>m_{n+1}]}\circ \Phi_n \\
= &N_{n+1} +\widetilde{R_{n+1}}. 
\end{aligned}
\eneq
Notice that by extracting all the terms of orders $m_n+1,\dots ,m_{n+1} $ from the equation above, one gets the cohomological equation \eqref{eq_homol}. 

\hfill {\bg $\Box$}

\bigskip

\subsection{Homological equation order by
  order.}\label{s_order_by_order}
Here we rewrite equation \eqref{eq_homol} as a (finite) set of
equations for each degree of $I$. Equations corresponding to degrees 
$m_n+1,\dots , m_{n+1}$ will formally determine $F_n$ (they are written out explicitly
in (\ref{homol-eq-by-order})). The rest of equations define $C_n$
(which is a part of the new remainder term).
Equating coefficients  with the same homogeneous degree in $I$
in both sides of \eqref{eq1} we obtain for the degrees from $m_n+1$ to
$m_{n+1}$
the following recursive formula (we write $m$ instead of $m_n$ for
typographic reasons): 
\beq\label{homol-eq-by-order}
\beal
&\{N_0, F^{[m]} \}+R^{[m+1]} =N^{ [m+1] }, \\
&\{N_0, F^{[m+1]} \} + \{N^{[3]}, F^{[m]} \}+ R^{[m+2]}= N^{[m+2]}, \\
&\{N_0, F^{[m+2]} \} + \{N^{[4]}, F^{[m]} \}+ \{N^{[3]}, F^{[m+1]} \}+ R^{[m+3]}= N^{[m+3]}, \\
&\dots \\
&\{N_0, F^{[2m-2]} \}+ \sum_{j=0}^{m-3}\{N^{[m-j]}, F^{[m+j]} \} + R^{[2m-1]}=N^{[2m-1]}.
\enal
\eneq
Recall that $2m_n-1=m_{n+1}$, see Sec. \ref{not_const}. From the formal
solvability 
we know that each of these equations has
a formal solution $F_n^{[m+j]}$. Of course, such a solution is not
unique. We will make the solution unique by prescribing the condition
$$
\int_{\T^d} F_n^{[m+j]}(I,\th) =0.
$$
As we will see, this normalization will allow us to get the estimates needed for the proof of the convergence.
The sum of the terms of orders $m_{n+1}+1, \dots , m_{n+1}+m_n-2$ (i.e., $2m_{n}, \dots , 3m_n-3$) that appear
in equation (\ref{homol-eq}) is denoted by $C_n$. In the notation $m=m_n$, we have: $C_n=C_n^{[2m, 3m-3]}$. 
The terms of the uniform degree satisfy 
\beq\label{C-by-order}
\beal
&C_n^{[2m]}=\{N_{[3]},  F^{[2m-2]} \}+ \{N_{[4]},  F^{[2m-3]} \}  +   \dots + \{N_{[m]},  F^{[m+1]} \}, \\
&C_n^{[2m+1]}=\{N_{[4]},  F^{[2m-2]} \}+   \{N_{[5]},  F^{[2m-3]} \}  +   \dots + \{N_{[m]},  F^{[m+2]} \}   \\
&\dots \\
&C_n^{[3m-3]}=\{N_{[m]},  F^{[2m-2]} \}.
\enal
\eneq
This can be written more compactly as
\beq\label{def_Cn}
C_n= \sum_{k=1}^{m-2} \{F^{[2m-1-k]},  \sum_{j=k+2}^{m}N^{[k+j]} \} . 
\eneq
This should be viewed as a definition of the remainder term
$C_n$.

\subsection{An important simplification.}\label{s_simplif}

In the case when the {\it normal form is an analytic function of 
$N_0(I)$} as in \eqref{N=B},
we have an important simplification. 
Denote 
\beq\label{notation_gj}
g_{2j}(I):= j b_{j} (N_0(I))^{j-1} \text{ and } 
g_{2j+1}(I)\equiv 0. 
\eneq 
Then for $j\in \N$ we have:
\beq\label{def_g_j}
\begin{aligned}
&\{N^{[2j]}, F\} = \{b_{j} (N_0)^{j}, F\}= j b_{j} (N_0)^{j-1} \{ N_0, F
\}=   
g_{2j}(I) \{ N_0, F \} \\
&\{N^{[2j+1]}, F\}= g_{2j+1}(I) \{ N_0, F \}\equiv 0.
\end{aligned}
\eneq

We formulate this as a lemma:
\blm\label{l_simplification}
If the normal form is an analytic function of 
$N_0(I)$ as in \eqref{N=B}, 
then equation \eqref{homol-eq-by-order} is equivalent to 
\beq\label{homol-eq-by-order_s}
\beal
&\{N_0, F^{[m]} \}+R^{[m+1]} =N^{ [m+1] }, \\
&\{N_0, F^{[m+1]} \} + g_{3}(I)\{N_0, F^{[m]} \}+ R^{[m+2]}= N^{[m+2]}, \\
&\{N_0, F^{[m+2]} \} + g_{4}(I) \{N_0, F^{[m]} \}+ g_{3}(I)\{N_0, F^{[m+1]} \}+ R^{[m+3]}= N^{[m+3]}, \\
&\dots \\
&\{N_0, F^{[2m-2]} \}+ \sum_{j=0}^{m-3}g_{m-j}(I)\{N_0, F^{[m+j]} \} + R^{[2m-1]}=N^{[2m-1]}.
\enal
\eneq
\elm
and 
\beq\label{def_Cn_G}
C_n= \sum_{k=1}^{m-2} \left( \{F^{[2m-1-k]}, \, N_0   \} \cdot \sum_{j=k+2}^{m} g_j \right). 
\eneq

\subsection{Homological equations  in majorants.} \label{s_majorants}
Here we study a simple recursive formula and estimate its terms. 
Later it will provide an important estimate of  $|\{N_0, F^{j}
\}|_{\rho_n,\rho_n}$. Here is the idea: suppose that in the Lemma
above for some
$\eps>0$, for all $j=0,\dots ,m$ we
have:
$$P_j:=|R^{[m+j]}|_{\rho_n,\rho_n} +
|N^{ [m+j] }|_{\rho_n,\rho_n}\leq \eps, \quad |g_{j}|_{\rho_n}\leq 1/4^j.$$
Define $S_j$ by the relations \eqref{homol-eq-by-order1} below. 
Then, by Lemma \ref{l_simplification}, for all $j=0,\dots ,m$ we have
$$
|\{N_0, F^{j} \}|_{\rho_n,\rho_n}\leq S_j.
$$

\blm\label{l_est_NF}
Given $ \eps >0$, 
suppose that for all $j=1, \dots ,m-1$ the numbers
$P_j$  satisfy
$$
0< P_j \leq \eps . 
$$
Let $S_j$ be defined recursively by equations 
\begin{equation}\label{homol-eq-by-order1}
\begin{aligned}
&S_{1} = P_{1}, \\
&S_{2} = P_{2} + \frac14 S_{1} , \\
&S_{3} = P_{3} + \frac14 S_{2}+ \frac1{4^2} S_{1} \\
&S_{4} = P_{4} + \frac14 S_{3}+ \frac1{4^2} S_{2}+ \frac1{4^3} S_{1} \\
&\dots \\
&S_{m-1} = P_{m-1} + \sum_{j=1}^{m-1} \frac1{4^j} S_{m-1-j} .\\
\end{aligned}
\end{equation}

Then for each  $j$ we have
$$
S_j \leq 2 \eps ,\quad j=1, \dots ,m-1.
$$
\elm
\bg \noindent{\it Proof.} \bk
By the formula for  $S^{[j]}$  above,
$$
S_j \leq P_j + \frac14 S_{j-1} + \frac14 ( S_{j-1}  -P_{j-1} ) = P_j+ 2 \frac14 S_{j-1} \leq P_j+  S_{j-1}/2 .
$$
This implies 
$$
S_j \leq  \sum_{k=0}^{j-1} 2^{-k} P_{j-k} \leq \eps \sum_{k=0}^{j-1}2^{-k}  < 2 \eps.
$$
\hfill {\bg $\Box$}

\section{Formal solution provides  analytic with estimates.}\label{s_one_step}
In this section we study a
 homological equation \eqref{homol-eq} below with an analytic right-hand side $Q(I,\th)$. 
 Assuming that it has a formal solution, we will find an analytic one,
 and estimate 
it in terms of the right hand side. Similar procedures appear in \cite{Ll}.

\blm\label{lemma-homol-eq0}
Let  $N_0(I)=I^{tr}\Om I$ where $\Om$ is a symmetric matrix with
$\det\Om\neq 0$, and let $Q(I,\th)$ be analytic in an annulus
$\A_{\rho,\sigma}$ for some $\rho$, $\sigma>0$. 
Suppose that the following equation  has a formal solution $\widetilde
F (I,\th)$:
\beq\label{homol-eq}
\{N_0,\widetilde  F \}= Q.
\eneq
Then equation \eqref{homol-eq} has an analytic  solution $F(I,\th)$, defined in
$\A_{\rho,\sigma}$,
and for any $0<\dt<\rho$, $0<\gamma<\sigma$ we have:
$$
|F |_{\rho -\dt, \sigma-\gamma} \leq c(d,\Om) \frac{1}{\dt \gamma^d} | Q |_{\rho, \sigma}  ,
$$
where $c(d,\Om)$ is a constant only depending on $d$ and $\Om$.

Moreover, if $Q(I,\th)$ is a homogeneous polynomial in $I$ with coefficients
depending on $\th$, then so is $F(I,\th)$. 

\elm

\bg \noindent{\it Proof. } \bk
Expanding $F$ formally into a Fourier series:
$F=\sum_{k\in \Z^d} \widehat F_k(I) e^{2\pi i\< k,\th\>}$, 
we get: 
$$
\{N_0,F\} = \sum_{j=1}^d F_{\th_j} (N_0)_{I_j} =  2\pi i \sum_{k\in \Z^d} \< k, 2\Om I\> \widehat F_k (I) e^{2\pi i \< k,\th\>}.
$$
Recall that $\Om$ is symmetric, so $\< k, \Om I\>=\< \Om k,  I\>$. Expressing $Q=\sum_{k\in \Z^d} \widehat Q_k(I) e^{2\pi i \< k,\th\>}$, 
we can rewrite equation (\ref{homol-eq}) as a series of equations indexed by $k$:
\beq\label{homol-eq-k}
\widehat Q_k(I)=4\pi i \< \Om k,  I\> \widehat F_k (I).
\eneq
If $\<k, \Om I\> \neq 0$, we can express $\widehat F_k = 
\widehat Q_k(I)/ (4\pi i \< \Om k,  I\>)$.

Since we have assumed existence of a formal solution of the
homological equation (\ref{homol-eq}) (and hence, a solution of (\ref{homol-eq-k}) for
each $k$), we have: 
$$
\< \Om k,  I\>=0 \Rightarrow \widehat Q_k(I)=0.
$$ 
Hence, for  $\< \Om k,  I\>=0$,
the equation is satisfied for any value of
$\widehat F_k(I)$. 
We define $\widehat F_k$ at these points by continuity. 
A way to do it is the following. 
Differentiate equation (\ref{homol-eq-k})
in the direction of $\Om k$:
$$
\< \Om k, \nabla \widehat Q_k(I) \>=   4\pi i \left( |\Om k |^2 \widehat F_k(I) + \< \Om k,  I\>
\< \Om k, \nabla \widehat F_k(I)\>  \right) , 
$$
where for a vector $v\in \R^d$ we denote $|v|^2=\sum_{j=1}^d v_j^2$.
For $\< \Om k,  I\> =0$, define $\widehat F_k(I)= \< \Om k , \nabla \widehat
Q_k(I)\> /(4\pi i |\Om k|^2)$. Summing up, we have defined a continuous function $\widehat F_k(I) $  by
$$
\widehat F_k (I) = \frac{1}{4\pi i} 
\begin{cases}
\< \Om k,  I\>^{-1} \widehat Q_k (I), \quad \< \Om k,  I\> \neq 0, \\
\frac{1}{|\Om k|^2}\< \Om k , \nabla \widehat Q_k(I)\> ,\ \  \< \Om k,  I\> =0.
\end{cases}
$$
Moreover,
since $\widehat F_k(I)$ is analytic in $\D_\rho \setminus \{\< \Om k,  I\> =0\}$ and
bounded in $\D_\rho$, it is analytic in $\D_{\rho}$.   
Notice that if in equation \eqref{homol-eq-k} $\widehat Q_k(I)$ is a homogeneous polynomial in $I$, 
then so is $\widehat F_k(I)$.

Now let us estimate the norm of the solution. Fix $0<\dt<\rho/2$, $0<\gamma<\sigma$. For each fixed 
$k\in\Z^d$,  we will estimate the corresponding $\widehat F_k(I)$ in two steps: first `$\delta/2$-close" to the resonant plane $\< \Om k,  I\> $, and then in the rest of $\D_{\rho -\dt} $. 

For the first step, let $\Pi_\dt= \{\<\Om k,I\>=0 \} \cap \D_{\rho-\dt}$ be the part of the 
resonant plane falling into $\D_{\rho-\dt}$. Notice that the orthogonal complement to this plane is formed 
by the vectors $\alpha e^{2 \pi i \phi } \Om k $, $\alpha\geq 0$, $\phi\in [0,1)$. Let 
$$
\Delta=\left\{ I= \alpha  \frac{\Om k }{|\Om k | } e^{2 \pi i \phi } \, \Big| \, \alpha < \dt  /2, \,\phi\in [0,1) \right\}
$$ 
be the complex disk 
of radius $\dt /2$ centered at zero and orthogonal to $\Pi_\dt$.  
Note that the restrictions of
$\widehat Q_k(I)$ and $\widehat F_k(I)$ to this disc 
are analytic.
Consider the ${\dt}/2$-neighbourhood $O_\dt$ of 
$\Pi_\dt$: 
$O_\dt=\bigcup_{I_0\in \Pi_\dt} (I_0+\Delta)$. Then and $O_\dt\subset \D_{\rho -\dt} $.

For each fixed $I\in O_\dt$ there exists $I_0 \in \Pi_\dt$ such that $I\in I_0+\Delta$.  We can estimate $|\widehat F_k (I)|$ by
the maximum modulus principle on the disk  $I_0+\Delta$. Namely, for $I$ lying on the boundary of this disk we have: $|\< \Om k,  I\>|= |\< \Om k,  I_0\>+ \<  \Om k, \dt \Om k/(2| \Om k|) \> |  = | \Om k|\dt/2$. Hence, for such $I$ we have
$$
|\widehat F_k (I)| \leq 
\frac{2 | \widehat Q_k |_\rho }{4\pi\dt |\Om k|}< \frac{ | \widehat Q_k |_\rho }{\dt |\Om k|}.
$$
As the second step in this estimate, consider $I\in \D_{\rho -\dt}\setminus O_\dt $. Here 
$|\< \Om k,  I\>| \geq  |\Om k| \dt / 2$, so $|\widehat F_k (I)| $ satisfies the same estimate as above.

By Cauchy estimates, we have: 
$$
|\widehat Q_k |_{\rho} \leq | Q |_{\rho,\sigma}e^{-|k|\sigma}.
$$
Since det$\, \Om\neq 0$,  there exists a constant $c(\Om)$ such that $|\Om k|\geq |k|/c(\Om)$ for all $k$. 
Then
$$
\begin{aligned}
|\widehat F_k |_{\rho-\dt} \leq \frac1{\dt |\Om k|} | \widehat Q_k |_{\rho}
\leq  c(\Om) \frac{e^{-\sigma |k|}} {\dt |k|}  | Q |_{\rho,\sigma}.
\end{aligned}
$$
Finally, for small $\dt$ and $\gamma$ we have:
$$
\begin{aligned}
|F |_{\rho-\dt,\sigma-\gamma}\leq &
\sum_{k\in \Z^d \setminus \{0\}} e^{(\sigma-\gamma) |k|} |\widehat F_k |_{\rho-\dt }
\leq  \frac{c(\Om)}{\dt }
\sum_{k\in\Z^d \setminus \{0\}}  \frac{e^{-\gamma |k|} }{|k|} | Q  |_{\rho,\sigma}\\
\leq &\frac{c(d,\Om)}{\dt \gamma^d } | Q |_{\rho,\sigma}  ,
\end{aligned}
$$
where $c(d,\Om)$ is a constant only depending on $d$ and $\Om$. The estimates above are very wasteful, but they are enough for our purposes.
\hfill {\bg $\Box$}

\section{Proof of Proposition \ref{lemma1}.}
Here we summarize the preparatory work to complete the proof of
Proposition \ref{lemma1}.
Let us return to the original problem.  
For a fixed $n$, let the necessary constants be as in Sec.\ref{not_const},
$|\widetilde{R_n} |_{\rho_n } \leq \dt_n^\kappa$, 
and
let $ g_{2j}(I)=j\, b_{j} \, (N_0(I))^{j-1}$ as in (\ref{notation_gj}).

\subsection{Estimate of  $|\{N_0,F_n \}|_{\rho_n, \rho_n}$ and 
$|C_n |_{\rho_n, \rho_n }$. }\label{s_est_NF}

For $j=1, \dots ,m_n-1$ denote 
$$
P_{j}:= |N^{[m_n+j]}|_{\rho_0 }      +|R^{[m_n+j]} |_{\rho_n } .
$$

By the choice of $\rho_0$, see Sec. \ref{not_const}, for all $j=1, \dots , m_n-1$  we have: 
$$
|g_{j}(I) |_{\rho_0} \leq  4^{-j},\quad   |N^{[m_n+j]}|_{\rho_0 } \leq  \dt_n^\kappa.
$$
Since for $j=1, \dots ,m_n-1$  we have
$|R^{[m_n+j]} |_{\rho_n } \leq |\widetilde{R_n} |_{\rho_n } \leq \dt_n^\kappa$,
and so for these values of $j$ we get
$$
P_j\leq 2\dt_n^\kappa.
$$
let $S_j$ be defined by (\ref{homol-eq-by-order1}). By Lemma \ref{l_est_NF}, 
for $j=1,\dots m-1$ we have $S_j \leq 2\eps $.
Equations (\ref{homol-eq-by-order_s}) imply that for $j=1,\dots m-1$ we have
\beq\label{eq22}
| \{N_0, F_n^{[m+j-1]} \} |_{\rho_n,\rho_n} \leq S_j\leq 2\eps =4\dt_n^\kappa.
\eneq
By linearity,
$$
| \{N_0, F_n \} |_{\rho_n,\rho_n} \leq  \sum_{j=1}^{m_n-1} | \{N_0, F_n^{[m_n+j-1]} \} |_{\rho_n,\rho_n}   
\leq   4m_n \dt_n^\kappa \leq  4\dt_n^{\kappa - 1}.
$$
The latter estimate follows from the definition of $m_n$ and $\dt_n$, see Section \ref{not_const}.

Moreover, by \eqref{def_Cn_G},
$$
\begin{aligned}
|C_n |_{\rho_n} &= \sum_{k=1}^{m-2} \left( S_{m-k} \,\sum_{j=k+2}^{m} G_j \right) \leq  \sum_{k=1}^{m-2} \left( S_{m-k} \,\sum_{j=k+2}^{\infty} 4^{-j} \right) \\
&\leq \frac{1}{3} \sum_{k=1}^{m-2} 4^{-(k+1)} S_{m-k} \leq \frac{1}{2} \eps=\dt_n^\kappa. 
\end{aligned}
$$
Hence,
\beq\label{estimate_Cn_S}
| C_n |_{\rho_n} \leq \dt_n^\kappa .
\eneq

\subsection{Estimates for $F_n$.}

Consider equation \eqref{eq22}.
Lemma \ref{lemma-homol-eq0} with $\rho=\sigma=\rho_n$,  $\dt=\gamma=\dt_n$ and $| Q |_{\rho, \sigma} \leq 4\dt_n^\kappa$, implies:
$$
| F_n^{[m+j-1]} |_{\rho_n-\dt_n, \rho_n-\dt_n}   \leq    4c(d,\Om) \dt_n^{\kappa - d-1} .
$$
Since $F_n=F_n^{[m_n, m_n+j-1]} $ where $m_n\leq \dt_n^{-1}$, we get: 
\begin{equation}\label{eq24}
| F_n |_{\rho_n - \dt_n , \rho_n-\dt_n} \leq \sum_{j=1}^{m_n-1} | F_n^{[m+j-1]} |_{\rho_n-\dt_n , \rho_n-\dt_n} \leq  m_n \, 4c(d,\Om) \dt_n^{\kappa - d-1}  
\leq  \dt_n^{\kappa - d-3} \leq  \dt_n^{3}.
\end{equation}
The latter estimate follows from the definition of $\kappa$, see Section \ref{not_const}.

\subsection{Estimates for $\Phi_n$.}
Here we prove that with $F_n$ as above, estimates (\ref{estPhi}) and (\ref{estPhi-1})
hold true. Indeed, 
the coordinate change  $\Phi_n = X_{F_n}^1$ is the time one map of the flow $X_{F_n}^t$
defined by the equations 
$$
\dot I = \partial_\th F_n (I,\th), \quad  \dot \th =-\partial_I F_n (I,\th).
$$
By \eqref{eq24}  and Cauchy estimates we get
\begin{equation}\label{eq_est_dF}
| \partial_I F_n |_{\rho_n - 2\dt_n , \rho_n-\dt_n} \leq  \dt_n^{2}, \quad  | \partial_\th F_n |_{\rho_n - \dt_n , \rho_n-2\dt_n} \leq  \dt_n^{2}.
\end{equation}
Then for any $t\leq 1$:
$$
 |X_{F_n}^t (I,\th)- (I,\th)|_{\rho_n - 3\dt_n , \rho_n- 3 \dt_n}       \leq  t \, \dt_n^{-1}  | F_n |_{\rho_n - 2\dt_n , \rho_n-2\dt_n} \leq \dt_n^{2}.
$$
\begin{equation}\label{eq_est_Phi}
X_{F_n}^t :\A_{\rho_n-3\dt_n,\rho_n-3\dt_n} \mapsto \A_{\rho_n-2\dt_n,\rho_n-2\dt_n} 
\end{equation}
In particular, since $ \Phi_n = X_{F_n}^1$, we get the desired formulas \eqref{estPhi} and \eqref{estPhi-1}.


\subsection{Estimate of the new remainder $\widetilde {R_{n+1}}$.}
\blm 
\label{lemma_est_Rn+1} For $F_n$ constructed above, estimate
(\ref{estRn}) holds:
$$
|\widetilde {R_{n+1}} |_{\rho_n - 3\dt_n , \rho_n-3\dt_n} <  4 \dt_{n}^\kappa.
$$
\elm

\bg \noindent{\it Proof.} \bk
By Lemma \ref{lem_homol_eq},  
$$
\widetilde {R_{n+1}}= A_n+B_n+C_n,
$$ 
where $A_n$, $B_n$ and $C_n$ are defined by
(\ref{not_AnBn}) and (\ref{not_Cn}).

{\bf Estimate of $A_n$: }
Using \eqref{eq_est_Phi},
we get:
$$
| \widetilde {R_{n}}^{[>m_{n+1}]}\circ  \Phi_n  |_{\rho_n-3\dt_n,\rho_n-3\dt_n} 
\leq |\widetilde {R_{n}} |_{\rho_n-2\rho_n,\rho_n-2\dt_n} \leq \dt_n^\kappa. 
$$

\medskip

{\bf Estimate of $C_n$:} We showed in section \ref{s_est_NF} that 
$$
|C_n|_{\rho_n,\rho_n}\leq \dt_{n}^\kappa.
$$

\medskip

{\bf Estimate of $B_n$:}
By \eqref{eq_est_dF},
$
| \partial_I F_n |_{\rho_n - 2\dt_n , \rho_n-\dt_n} \leq  \dt_n^{2}$ and $  | \partial_\th F_n |_{\rho_n - \dt_n , \rho_n-2\dt_n} \leq  \dt_n^{2}.
$
By (\ref{estRn}) 
$$
|R_{n} |_{\rho_n,\rho_n} \leq |\widetilde {R_{n}} |_{\rho_n, \rho_n} \leq \dt_n^\kappa. 
$$
This implies, using Cauchy estimates, that
$$
| \{ R_n,F_n\} |_{\rho_n - 2\dt_n , \rho_n-2\dt_n} \leq \dt_n^\kappa.
$$

Notice that,  by formulas \eqref{eq_homol}   and \eqref{not_Cn}, we have
$ \{ N_n, F_n \} = R_{n} +N_{n} - N_{n-1} +C_n$.

By \eqref{est_N},
$$
|N_{n} - N_{n-1}|_{\rho_0,\rho_0} =\sum_{j=1}^{m_n}N^{[m_n+j]}\leq m_n \dt_n^{\kappa+1}\leq \dt_n^{\kappa}.
$$
and therefore
$$
| \{ N_n, F_n \} |_{\rho_n,\rho_n}  = |R_{n} |_{\rho_n,\rho_n} + |N_{n} - N_{n-1} |_{\rho_n,\rho_n}+|C_n|_{\rho_n,\rho_n}
\leq 3 \dt_{n}^\kappa.
$$
Combining the above estimates, we get
$$
| \{ \{ N_n, F_n \} ,F_n\} |_{\rho_n - 2\dt_n , \rho_n-2\dt_n} \leq \dt_n^\kappa,
$$
Since, by \eqref{eq_est_Phi}, for any $t\leq 1$ we have $X_{F_n}^t :\A_{\rho_n-3\dt_n,\rho_n-3\dt_n} \mapsto \A_{\rho_n-2\dt_n,\rho_n-2\dt_n} $,
we obtain
$$
| \{ \{ N_n, F_n \} +R_n,F_n\}  \circ X_{F_n}^t |_{\rho_n-3\dt_n,\rho_n-3\dt_n} \leq 
| \{ \{ N_n, F_n \} +R_n,F_n\}  |_{\rho_n-2\dt_n,\rho_n-2\dt_n}     \leq 2 \dt_n^\kappa.
$$

\hfill {\bg $\Box$}

\bigskip

Here we get the desired estimate for the remainder term. We have proved above that 
$$
|\widetilde {R_{n+1}} |_{\rho_n-3\dt_n,\rho_n-3\dt_n}< 4 \dt_{n}^{\kp}
$$
Recall that 
$\widetilde {R_{n+1}}=\widetilde {R_{n+1}}^{[>m_{n+1}]}$.
By Lemma \ref{lemma_est_Rn+1_final}  proved below, this implies the desired estimate
$$
|\widetilde {R_{n+1}} |_{\rho_{n+1},\rho_{n+1}}<  \dt_{n+1}^{\kp}
$$
This finishes the proof of Proposition \ref{lemma1}, and hence Theorem 1 (as explained
in the introduction). \hfill {\bg $\Box$}
\medskip

\blm 
\label{lemma_est_Rn+1_final} 
Suppose that the constants $\kp$, $b$, $\dt_n$, $q_n$, $\rho_n$ are defined in Section \ref{not_const}, and
an analytic function $G(I,\theta)$ satisfies  $G=G^{[>m_{n+1}]}$, and
$$
|G |_{\rho_n-3\dt_n,\rho_n-3\dt_n}< 4 \dt_{n}^{\kp}.
$$
Then 
$$
|G|_{\rho_{n+1},\rho_{n+1}}<  \dt_{n+1}^{\kp}.
$$
\elm

\bg \noindent{\it Proof.  } \bk 
By  the definition of $\kp$ in Section \ref{not_const} we have:
$q_n^{m_{n+1}+1}=q_{n}^{2^{n+1}+2} < q_{n}^{2^{n+1}} = 2b= 2^{-\kp-2}$. 
Also recall that $\dt_{n+1}=2^{-1}\dt_{n}$.

Since $G$ starts with terms of degree
$m_{n+1}=2^{n+1}+2$, we have:
$$
|G |_{q_n(\rho_n-3\dt_n),q_n(\rho_n-3\dt_n)}< q_n^{2^{n+1}+2} \,4 \dt_{n}^{\kp}
\leq  2^{-\kp-2} \, 4 \dt_n^{\kp}   \leq  \dt_{n+1}^{\kp}. 
$$
\hfill {\bg $\Box$}

\Addresses

\end{document}